\documentclass[reqno,11pt]{amsproc}
\newtheorem{theo}{Theorem}

\newtheorem{lem}{Lemma}

\newtheorem{df}{Definition}
\newtheorem{hyp}{Hypothesis}
\newcommand\eps\varepsilon
\newcommand\ph\varphi
\newcommand\kap\varkappa
\newcommand\R {\mathbb{R}}

\begin{document}

\title[Abstract version of the Cauchy-Kowalewski Problem]
{On some abstract version of the Cauchy-Kowalewski Problem}

\author[Oleg Zubelevich]{Oleg Zubelevich\\ \\
\em Department of Differential Equations\\
Moscow State Aviation Institute\\
Volokolamskoe Shosse 4, 125871, Moscow, Russia\\
\rm E-mail: ozubel@yandex.ru}

\address{Department of Differential Equations
Moscow State Aviation Institute
Volokolamskoe Shosse 4, 125871, Moscow, Russia}
\email{ozubel@yandex.ru}
\curraddr{2-nd  Krestovskii Pereulok 12-179, 129110, Moscow, Russia}

\thanks{Partially supported by grants RFBR 02-01-00400, INTAS 00-221.}

\begin{abstract}We consider an abstract version of the
Cauchy-Kowalewski Problem with the right hand side being
free from the Lipschitz type conditions and prove the existence theorem.
\end{abstract}

\maketitle
\numberwithin{equation}{section}
\newtheorem{theorem}{Theorem}[section]
\newtheorem{lemma}[theorem]{Lemma}
\newtheorem{definition}{Definition}[section]

\section{Introduction}
There are two most standard existence theorems
in the theory of ODE: the Cauchy-Picard existence and uniqueness
theorem and the Peano existence theorem. The Cauchy-Picard theorem states
that if the right hand side of ODE satisfies the Lipschitz conditions
then initial value problem has unique solution. The proof of this theorem
is based on the contraction mapping principle. The Peano theorem states
that for existence of a solution it is sufficient to have only
continuousness of the right hand side.
This theorem is proved by means of
compactness considerations with the help of the Arzela-Ascoli theorem.

The case of initial value problem for PDE in the abstract setup has been
studied by many authors and there are existence and uniqueness theorems
proved under the assumptions of Lipschitz type conditions.

An abstract form of the Cauchy-Kowalewski Problem was first considered by
T. Yamanaka in \cite{Yamanaka15} and L. Ovsjannikov  \cite{Ovsjannikov10} in
the linear case. Some another aspects of the linear Cauchy-Kowalewski Problem was
exposed by J. Treves \cite{Treves12}.

In \cite{Nirenberg8} L. Nirenberg obtained the existence and uniqueness
theorem for  the abstract nonlinear Cauchy-Kowalewski Problem.
The proof of Nirenberg's theorem uses an iteration procedure of Newtonian
type and based on ideas of the KAM theory. In Nirenberg's theorem  it is
assumed that the right hand side of the problem is a strong differentiable mapping.

T. Nishida in \cite{Nishida} simplified the iteration procedure
and stated that in Nirenberg's theorem it is possible
to replace strong differentiability with the Lipschitz type conditions.

In \cite{Safonov} M. Safonov
gave a proof of Nishida's theorem by constructing a suitable Banach space of functions
and then using the contraction mapping principle.

In present paper we consider a topological aspect of the abstract
nonlinear Cauchy-Kowalewski Problem and prove the Peano type existence theorem.

We assume that the right hand side of the equation
depends on two arguments: it is bounded and
continuous in the first argument (pure Peano's case) and convex in the
second one. Such a setup includes  quasilinear PDE as a special case. This theorem
is not deduced from Nishida's result or quasilinear
versions of the Cauchy-Kowalewski Problem since  the Lipschitz type conditions
are not applied.

The main tools we use is Browder's generalization of the Schauder fixed
point theorem and a topological construction close to Safonov's one.

\section{Main theorem}
Let $\{(E_s,\|\cdot\|_s)\}_{0<s<1}$ be a scale of Banach spaces:
\begin{equation}
\label{scale}
 E_{s+\delta}\subseteq E_s,\quad
\|\cdot\|_s\le\|\cdot\|_{s+\delta},\quad s+\delta<1,\quad
\delta>0.\end{equation}
We assume that all embeddings (\ref{scale}) are compact.
Such an assumption always holds for the scales of analytic functions.

Let $B_s(r)=\{u\in E_s\mid \|u\|_s<r\}$ be an open ball of $E_s$ and let
$\overline{B}_s(r)$ be its closure.

The main object of our study is the following Cauchy-Kowalewski problem:
\begin{equation}
\label{problem}
u_t=A(t,u,u)+h(t,u),\quad u\mid_{t=0}=0.
\end{equation}
For some positive constant $T>0$ the mappings
$$A:[0,T]\times E_{s+\delta}\times E_{s+\delta}\to E_s,\quad h:[0,T]
\times E_{s+\delta}\to E_{s},\quad \delta>0,\quad s+\delta<1$$
are continuous and there are positive constants $R,M,K$
such that if $u,v\in \overline{B}_{s+\delta}(R)$ then the inequalities hold:
\begin{equation}
\label{as_1}
\|A(t,u,v)\|_s\le \frac{M\|v\|_{s+\delta}}{\delta},\quad
\|h(t,u)\|_s\le K,\quad \delta>0,\quad s+\delta<1.
\end{equation}
Let the mapping $A$ be convex in the third argument: for all
$u,v,w\in B_{s+\delta}(R)$ and $0\le\lambda\le 1$
we have
\begin{equation}
\label{as_conv}
\|A(t,u,\lambda v+(1-\lambda)w)\|_s\le \lambda\|A(t,u,v)\|_s+(1-\lambda)\|A(t,u,w)\|_s.
\end{equation}
For example, if the mapping $A$ is linear in the third argument then the above inequality holds.

\begin{theo}
\label{main_th}There exists such a large constant $a>0$
that problem (\ref{problem}) has a solution
$$u(t)\in \bigcap_{1-s-\tau a>0}C([0,\tau],\overline{B}_{s}(R)).$$
\end{theo}
\begin{hyp}The condition of convexness (\ref{as_conv}) is necessary: there
must be an example of such  mappings $A$ and $h$ that satisfy all the above conditions
except (\ref{as_conv}) and problem (\ref{problem}) does not have the
solution.
\end{hyp}

This theorem does not reduce to the Nishida  result \cite{Nishida}. Nishida's theorem
uses some kind of the Lipschitz condition:
$$
\|f(t,u')-f(t,u'')\|_s\le \frac{M}{\delta}\|u'-u''\|_{s+\delta},\quad
u',u''\in B_{s+\delta}(R),\quad s+\delta<1,$$where $f$ is a right hand side of the problem.

In the case under consideration we separate the arguments of the mapping $A$.
It is bounded in the second argument  and unbounded in the third one.
Thus it is sufficient to have only continuousness in the second argument and
linearity or convexness in the third one.

If the mapping $A$ equals to zero
identically then Theorem \ref{main_th} is a direct generalization from the
finite dimensional case to the scale of Banach spaces of the Peano
existence theorem.

There is no reason to expect uniqueness in Theorem \ref{main_th}: even in
the case of ordinary differential equations there are systems
with continuous (but not Lipschitz) right-hand side that do not have the
uniqueness.

Before starting to prove Theorem \ref{main_th} we must build some

\section{Preliminary topological construction}Introduce a triangle:
$$\Delta=\{(\tau,s)\in\R^2\mid \tau>0,\quad 0<s<1,\quad 1-s-\tau a>0\}.$$
Consider a seminormed space $E=\bigcap_{(\tau,s)\in\Delta}C([0,\tau],E_{s})$
with a family of  norms:
$$\|u\|_{\tau,s}=\max_{0\le t\le \tau}\|u(t)\|_s.$$
Obviously, these norms satisfy the following inequalities:
\begin{equation}
\label{norm_estmates}
\|\cdot\|_{\tau,s}\le\|\cdot\|_{\tau+\delta,s},\quad
\|\cdot\|_{\tau,s}\le\|\cdot\|_{\tau,s+\delta},\quad
\delta>0.\end{equation}

The space $E$ is a topological space with a basis of the topology given by
the open balls:
$$B_{\tau,s}(R)=\{u\in E\mid \|u\|_{\tau,s}<R\}.$$
\begin{df}
A set $G\subseteq E$ is said to be uniformly continuous if for all $\eps>0$ and
for all $(\tau,s)\in \Delta$  there is
$\delta=\delta(\eps,\tau,s)>0$ such that if $t_1,t_2\in [0,\tau]$ and $|t_1-t_2|<\delta$ then
$$\sup_{u\in G}\|u(t_1)-u(t_2)\|_s< \eps.$$

A set $G\subseteq E$ is said to be bounded if there are such constants
$M_{\tau,s}$ that for all $u\in G$ we have $\|u\|_{\tau,s}\le M_{\tau,s}$.
\end{df}

Recall the Arzela-Ascoli lemma \cite{Schwartz}:
\begin{lem}
\label{Ascoli}
Let $H\subset C([0,T],X)$ be
a set in the space of continuous functions
with values in a Banach space $X$.
Assume that the set $H$ is closed, bounded, uniformly continuous and
for every $t\in [0,T]$ the set $\{u(t)\in X\}$ is compact in the space
$X$. Then the set $H$ is compact in the space $C([0,T],X)$.
\end{lem}
Obviously there is a similar compactness criterium for the
space $E$.
\begin{lem}
\label{Ascoli_E}
If a closed set $G\subseteq E$ is uniformly continuous and bounded then it is
compact.\end{lem}
\begin{proof}Let $(\tau,s)$ be an arbitrary point of $\Delta$.
Since the set $G$ is bounded and uniformly continuous in the space
$C([0,\tau],E_{s+\delta})$, by the Lemma \ref{Ascoli}
it is compact in the space $C([0,\tau],E_{s})$.
Thus every sequence
$\{u_k\}_{k\in \mathbb{N}}\subset G$ contains a subsequence that converges with
respect to the norm $\|\cdot\|_{\tau,s}$. So it remains to prove
that there is a subsequence of $\{u_k\}$ that converges by all the
norms $\|\cdot\|_{\tau,s}$ at once.

Consider a set $\Delta_\mathbb{Q}=\Delta\bigcap\mathbb{Q}^2$. This set is
countable and let $\gamma:\mathbb{N}\to \Delta_\mathbb{Q}$ be a
corresponding bijection.

Let $\{u_k^1\}\subseteq \{u_k\}$ be a subsequence that converges by the
norm $\|\cdot\|_{\gamma(1)}$. By the above arguments there is a subsequence
$\{u_k^2\}\subseteq \{u_k^1\}$ that converges by the norm $\|\cdot\|_{\gamma(2)}$ etc.

The diagonal sequence $\{u_k^k\}$ converges by  the norms
$\{\|\cdot\|_{\gamma(k)}\}_{k\in \mathbb{N}}$. Then due to inequalities
(\ref{norm_estmates}) it converges in all the norms.

\end{proof}

In the conclusion we formulate a generalized version of the Schauder fixed
point theorem.
\begin{theo}[\cite{Browder}]\label{Schau-Tich}Let $W$ be a compact and convex subset of the
seminormed space $E$. Then a continuous mapping $f:W\to W$ has a fixed
point $\hat u$ i.e. $f(\hat u)=\hat u$.
\end{theo}

\section{Proof of Theorem \protect{\ref{main_th}}}
Problem (\ref{problem}) is obviously equivalent to the following one:
$$u(t)=F(u)=\int_0^tA(\xi,u(\xi),u(\xi))+h(\xi,u(\xi))\,d\xi.$$
So, we seek for a fixed point of the mapping $F$.

Let $S\subset E$ be a set that consists of such elements $v$ that
satisfy the following conditions:
\begin{equation}
\label{s1}\|v(t)\|_s\le R,\end{equation}
for all $u\in E$ such that $\|u(t)\|_s\le R$ we have
\begin{equation}\label{s2}
\|A(t,u(t),v(t))\|_s\le\frac{1}{\sqrt{1-at-s}},\end{equation}
for $(t_i,s)\in\Delta,\quad i=1,2$ we have
\begin{equation}\label{s3}
\|v(t_1)-v(t_2)\|_s\le
\Big(K+\frac{2}{a}\Big)|t_1-t_2|.\end{equation}
Note that the set $S$ is nonvoid: $0\in S$.

The set $S$ is convex by inequality (\ref{as_conv}) and it is compact by
Lemma \ref{Ascoli_E}.

Thus if we show that
\begin{equation}
\label{incl}
F(S)\subseteq S
\end{equation}
then the Proof will be conclude by
applying Theorem \ref{Schau-Tich} to the mapping $F$ and the set $S$.

Inclusion (\ref{incl}) is developed in the following order: first observing that
$$t<\frac{1}{a}$$
 we
verify that the mapping $F$ preserves inequality (\ref{s1}) then it implies
the same regarding (\ref{s2}) and then (\ref{s3}). Every step we choose
constant $a$ to be sufficiently large.

To illustrate this technique we assume that the
preserving of inequality (\ref{s1}) by the mapping $F$ has already been
checked up and
verify inequality (\ref{s2}). Other conditions are easier to obtain and they
go in the same manner.

So let $v\in S$. Evaluate by formulas (\ref{as_1}):

\begin{align}
\|A(t,u,F(v))\|_s &\le
M\int_0^t\frac{\|A(\xi,v(\xi),v(\xi))\|_{s+\delta}+\|h(\xi,v(\xi))\|_{s+\delta}}
{\delta}\,d\xi\nonumber\\
&\le M\int_0^t\frac{1}{\delta\sqrt{1-a\xi-s-\delta}}+\frac{K}{\delta}\,d\xi.\label{kn}
\end{align}

Substituting to the last formula $\delta=(1-a\xi-s)/2$
 we obtain that the value of expression (\ref{kn}) does not
exceed
$$
\frac{M}{a\sqrt{1-at-s}}\Big(2^{\frac{5}{2}}+K\Big).$$
So for $a$ greater or equals to
$M(2^{\frac{5}{2}}+K)$
 the mapping $F$ preserves inequality
(\ref{s2}).

Theorem \ref{main_th} is proved.

 \end{document}